\documentclass[a4paper,oneside,english]{amsart}
\usepackage[T1]{fontenc}
\usepackage[latin9]{inputenc}
\usepackage{textcomp}
\usepackage{amsthm}
\usepackage{graphicx}
\usepackage{amssymb}

\makeatletter
\numberwithin{equation}{section} 
\numberwithin{figure}{section} 
\theoremstyle{plain}
\theoremstyle{plain}
\newtheorem{thm}{Theorem}
  \theoremstyle{plain}
  \newtheorem{prop}[thm]{Proposition}
  \theoremstyle{definition}
  \newtheorem{defn}[thm]{Definition}
  \theoremstyle{remark}
  \newtheorem{rem}[thm]{Remark}
  \theoremstyle{definition}
  \newtheorem{condition}[thm]{Condition}
  \theoremstyle{definition}
  \newtheorem*{condition*}{Condition}
  \theoremstyle{plain}
  \newtheorem{cor}[thm]{Corollary}
  \theoremstyle{plain}
  \newtheorem{lem}[thm]{Lemma}
  \theoremstyle{remark}
  \newtheorem*{acknowledgement*}{Acknowledgement}

\DeclareMathOperator{\bd}{Bd}

\@ifundefined{showcaptionsetup}{}{%
 \PassOptionsToPackage{caption=false}{subfig}}
\usepackage{subfig}
\makeatother

\usepackage{babel}

\begin{document}

\title{Equivariant closure operators and trisp closure maps}

\author{Juliane Lehmann}

\address{Fachbereich Mathematik, Universität Bremen, 28359 Bremen, Germany}

\email{jlehmann@math.uni-bremen.de}
\begin{abstract}
A trisp closure map $\phi$ is a special map on the vertices of a
trisp $T$ with the property that $T$ collapses onto the subtrisp
induced by the image of $\phi$. We study the interaction between
trisp closure maps and group operations on the trisp, and give conditions
such that the quotient map is again a trisp closure map. Special attention
is on the case that the trisp is the nerve of an acyclic category,
and the relationship between trisp closure maps and closure operators
on posets is studied.
\end{abstract}



\maketitle

\section{Introduction}

This paper aims to bring together the subjects of two different papers:
\cite{BK05}, where Babson and Kozlov studied quotients of posets
and conditions under which these commute with the nerve functor; and
\cite{Ko08b}, where Kozlov introduced trisp closure maps, which are
a compact certificate for the collapsibility of a trisp to a certain
subtrisp.

Section~\ref{sec:Definitions} gives definitions for the needed objects:
acyclic categories and their nerves, quotients of trisps, quotients
of acyclic categories. Trisp closure maps get introduced in Section~\ref{sec:Cl-maps-on-trisps},
and there is a discussion on the relationship between them and closure
operators on posets. The main results of this paper can be found in
Section~\ref{sec:Cl-maps-on-quot-trisps}: Conditions that are sufficient
to guarantee the regularity of a quotient trisp also guarantee that
the quotient of a trisp closure map is a closure map on the quotient
trisp. In particular, these conditions are always fulfilled if the
trisp is the nerve of an acyclic category, with the group operation
induced by an operation on the acyclic category. Conversely, a trisp
closure map on a quotient trisp can be lifted whenever the lifted
map can be defined in a natural way and the original trisp is actually
a simplicial complex. In Section~\ref{sec:cl-maps-on-quot-poset},
we shed more light on the connection between closure operators and
trisp closure maps. We now consider the quotient of a poset, equipped
with a closure operator. As this quotient is taken in the category
of acyclic categories, we cannot expect the quotient of the closure
operator to be a closure operator again. But we do still obtain a
trisp closure map on the nerve of the poset quotient.

Finally, we apply our results in Section~\ref{sec:Applications}
to the $\mathcal{S}_{n}$-action on the barycentric subdivision of
the complex of disconnected graphs introduced by Vassiliev.

\section{\label{sec:Definitions}Definitions}

\subsection{Acyclic categories and posets}

In this paper, let $C$ be a finite acyclic category, that is, a finite
category where only identity morphisms have inverses. That means that
$C$ can be pictured with all arrows pointing upward; we write $s(m)$
for the source object of $m\in\mathcal{M}(C)$, $t(m)$ for the target
object. We write $x\parallel y$, if $\mathcal{M}(x,y)=\mathcal{M}(y,x)=\emptyset$
for two objects $x,y$. By $P$ we always denote a finite poset, which
we understand here as an acyclic category where for any two objects
$x,y$ there is at most one morphism $x\rightarrow y$, written $x\leq y$.
An AC-map on $C$ is a functor $\phi:C\rightarrow C$; if $C$ is
actually a poset, then this notion coincides with that of an order-preserving
map.

Let $G$ be a group acting on $C$. An action of $G$ on $C$ is a
functor $\mathcal{A}_{G}$ from the category $G$ to the category
of acyclic categories. This means that each group element $g$ acts
as an automorphism $\mathcal{A}_{g}$ on $C$; we write $ga$ instead
of $\mathcal{A}_{g}a$ for some morphism or object $a$. In particular,
such a group action here is always horizontal, meaning that $gx\neq x$
implies $gx\parallel x$ for all objects $x$. An AC-map $\phi$ on
$C$ is $G$-equivariant if it commutes with $\mathcal{A}_{g}$ for
all group elements $g$.

\subsection{Nerves of acyclic categories}

For the definition of a trisp, see for example \cite[Ch. 2]{Ko08a};
basically it is a generalization of an abstract simplicial complex,
where on the one hand there may be multiple simplices with the same
vertex set, and on the other hand there is an order on the simplices
compatible with taking boundaries. A trisp is regular if the number
of distinct vertices of each simplex equals the dimension of the simplex
plus 1. To each acyclic category $C$, we can associate a trisp $\Delta(C)$,
called the nerve of $C$. The 0-simplices (or vertices) of $\Delta(C)$
are the objects of $C$; the $t$-simplices are chains of $t$ composable
non-identity morphisms, e.g. $\sigma=a_{0}\stackrel{m_{1}}{\rightarrow}a_{1}\stackrel{m_{2}}{\rightarrow}a_{2}\rightarrow\ldots\stackrel{m_{t}}{\rightarrow}a_{t}$.
The boundary simplices of $\sigma$ are as follows: $\partial_{0}\sigma=a_{1}\stackrel{m_{2}}{\rightarrow}a_{2}\rightarrow\ldots\stackrel{m_{t}}{\rightarrow}a_{t}$,
$\partial_{i}\sigma=a_{0}\stackrel{m_{1}}{\rightarrow}\ldots\stackrel{m_{i-1}}{\rightarrow}a_{i-1}\stackrel{m_{i+1}\circ m_{i}}{\rightarrow}a_{i+1}\rightarrow\ldots\stackrel{m_{t}}{\rightarrow}a_{t}$,
$\partial_{t}\sigma=a_{0}\stackrel{m_{1}}{\rightarrow}a_{1}\stackrel{m_{2}}{\rightarrow}a_{2}\rightarrow\ldots\stackrel{m_{t-1}}{\rightarrow}a_{t-1}$.
Thus, the minimal vertex of $\sigma$ is $a_{0}$. The nerve of any
acyclic category is a regular trisp, so in the following we will be
concerned only with regular trisps without explicit mention.

Even more specially, nerves of acyclic categories are flag complexes.
That is, they are maximal under the condition that 1-skeleta of simplices
are unique. Thus, the trisp as a whole is uniquely determined by its
1-skeleton.

There is an induced action of $G$ on $\Delta(C)$, by \cite[Prop. 14.4]{Ko08a}
$\Delta(C)/G$ is again a regular trisp and the simplices are exactly
the orbits of simplices of $\Delta(C)$, e.g. $G\sigma=G(a_{0}\stackrel{m_{1}}{\rightarrow}a_{1}\stackrel{m_{2}}{\rightarrow}a_{2}\rightarrow\ldots\stackrel{m_{t}}{\rightarrow}a_{t})$
and the boundary maps work as expected: $\partial_{i}(G\sigma)=G(\partial_{i}\sigma)$.

\subsection{\label{sub:The-quotient-C/G}The quotient $C/G$}

Following \cite{BK05}, we define the quotient $C/G$ of $C$ by its
$G$-action as the colimit of the functor $\mathcal{A}_{G}$. In our
finite case here one can give an explicit description of $C/G$. The
objects $[a]$ are simply the $G$-orbits of objects $a$ of $C$.
The morphisms $[x]$ are equivalence classes of morphisms of $C$,
where the equivalence relation is induced by the $G$-action and composition,
with the transitive closure taken. Stated precisely, we have $[x]=[y]$
for $x,y\in\mathcal{M}(C)$ if there exist $z_{1},z_{2},\ldots,z_{n}\in\mathcal{M}(C)$,
$z_{1}=x,z_{n}=y$, and decompositions $z_{i}=z_{i,t_{i}}^{+}\circ z_{i,t_{i}-1}^{+}\circ\ldots\circ z_{i,1}^{+}$
for $i=1,2,\ldots,n-1$ and $z_{i}=z_{i,t_{i-1}}^{-}\circ z_{i,t_{i-1}-1}^{-}\circ\ldots\circ z_{i,1}^{-}$
for $i=2,3,\ldots,n$, such that $Gz_{i,j}^{+}=Gz_{i+1,j}^{-}$ for
all appropriate $i,j$.

As $C/G$ is again an acyclic category, we can consider its nerve
$\Delta(C/G)$. The universal property of colimits guarantees the
existence of a canonical map $\lambda:\Delta(C)/G\rightarrow\Delta(C/G)$,
with $\lambda(G(a_{0}\stackrel{m_{1}}{\rightarrow}a_{1}\stackrel{m_{2}}{\rightarrow}a_{2}\rightarrow\ldots\stackrel{m_{t}}{\rightarrow}a_{t}))=([a_{0}]\stackrel{[m_{1}]}{\rightarrow}[a_{1}]\stackrel{[m_{2}]}{\rightarrow}[a_{2}]\rightarrow\ldots\stackrel{[m_{t}]}{\rightarrow}[a_{t}])$.
On the 0-skeleta, $\lambda$ is an isomorphism, as $[a]=Ga$ are the
vertices of $\Delta(C/G)$ and $\Delta(C)/G$, respectively. Necessary
and sufficient conditions for this map to be an isomorphism have been
studied in \cite{BK05}; in particular $\lambda$ is always surjective.
As we will make heavy use of this fact, the proof is repeated here.
\begin{prop}
\cite[Prop. 3.1]{BK05} Let $C$ be an acyclic category with a $G$-action.
Then the canonical map $\lambda:\Delta(C)/G\rightarrow\Delta(C/G)$
is surjective.\end{prop}
\begin{proof}
Let $[a_{0}]\stackrel{[m_{1}]}{\rightarrow}[a_{1}]\stackrel{[m_{2}]}{\rightarrow}[a_{2}]\rightarrow\ldots\stackrel{[m_{t}]}{\rightarrow}[a_{t}]$
be a simplex of $\Delta(C/G)$. We will construct a $\lambda$-preimage
inductively. For $t=0$ the only possible choice as mentioned above
is $Ga_{0}$. If we have found $b_{0}\stackrel{n_{1}}{\rightarrow}b_{1}\stackrel{n_{2}}{\rightarrow}\ldots\stackrel{n_{t-1}}{\rightarrow}b_{t-1}$
with $[n_{i}]=[m_{i}]$ for all $i=1,\ldots,t-1$, implying $b_{t-1}\in[a_{t-1}]$,
there exists $g\in G$ such that $b_{t-1}=g\cdot s(m_{t})$. Thus
choosing $n_{t}=g\cdot m_{t}$ yields an extension to the composable
morphism chain, with $t(n_{t})\in t([m_{t}])=[a_{t}]$.
\end{proof}

\section{\label{sec:Cl-maps-on-trisps}Trisp closure maps}
\begin{defn}
A \emph{closure operator} on a poset $P$ is an order-preserving map
$\phi:P\rightarrow P$ with $\phi^{2}=\phi$. It is a descending (ascending)
closure operator, if $\phi(x)\leq x$ ($\phi(x)\geq x$) for all $x\in P$. 
\end{defn}
It is a well-known fact that a monotone closure operator $\phi$ on
$P$ induces a strong deformation retract from the order complex of
$P$ to the order complex of $\phi(P)$; see e.g. \cite[Cor. 10.12]{Bj96}.
Even more, $\Delta(P)$ actually collapses onto $\Delta(\phi(P))$
(see \cite[Th. 2.1]{Ko04}). Forgetting about the underlying poset
and just considering a trisp, we arrive at the central definition
of this paper, due to Kozlov (\cite{Ko08b}):
\begin{defn}
A \emph{trisp closure map} on a trisp $T$ is a partition of the vertex
set of $T$ into the \emph{blue} vertices $B$ and the \emph{red}
vertices $R$, together with a map $\phi:B\rightarrow R$ with the
following property: Let $\sigma$ be a simplex of $T$ containing
at least one blue vertex; let $b$ be the minimal blue vertex of $\sigma$.
Then either

\begin{itemize}\item$\phi(b)$ is a vertex of $\sigma$, and removing
$\phi(b)$ yields another (unique) simplex of $T$

\item or $\phi(b)$ is not a vertex of $\sigma$, then there exists
a unique vertex $\tau$ of $T$ that contains $\phi(b)$ as a vertex
and $\sigma$ as a boundary simplex of codimension 1.\end{itemize}
\end{defn}
Replacing {}``minimal'' with {}``maximal'' yields no conceptual
difference, and all statements in this paper which do not explicitely
mention the choice made still hold.

This definition is made worthwile by the following theorem due to
Kozlov.
\begin{thm}
\cite[Thm. 2.2]{Ko08b} Let $T$ be a regular trisp with a trisp closure
map $\phi:B\rightarrow R$. Then $T$ collapses on the subtrisp $T_{R}$,
consisting of those simplices of $T$ that contain only red vertices.\end{thm}
\begin{rem}
\label{rem:desc-cl-map->tr-cl-map}In particular, any descending (ascending)
closure operator $\phi$ on a poset $P$ induces a trisp closure map
$\bar{\phi}$ on $\Delta(P)$ with minimal (maximal) vertices chosen,
by setting $R=\phi(P),B=P\backslash R,\bar{\phi}=\phi_{|B}$, implying
that $\Delta(P)$ collapses onto $\Delta(\phi(P))$ (\cite[Cor. 2.5]{Ko08b}).
If $P$ carries a $G$-action and $\phi$ is $G$-equivariant, then
$R$ and $B$ are closed under $G$, since with $r=\phi(p)\in R$
also $gr=g\phi(p)=\phi(gp)$ is in $R$ for all group elements $g$.
\end{rem}
\begin{figure}

\hfill{}\subfloat[$\phi:b_{1}\mapsto r_{1},b_{2}\mapsto r_{2},b_{3}\mapsto r_{2},b_{4}\mapsto r_{3}$]{

\includegraphics{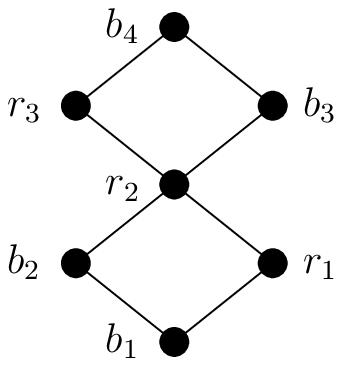}}\hfill{}\subfloat[$\bar{\phi}:b\mapsto r_{2}$]{

\includegraphics{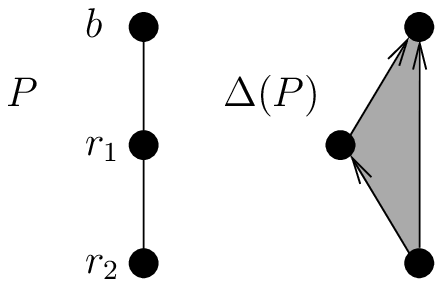}}\hfill{}

\caption{\label{fig:counterex}}

\end{figure}

The relationship between trisp closure maps and poset closure operators
is not as close as it appears on the first glance though. For example,
Figure~\ref{fig:counterex}a gives a closure operator on a poset
that does not induce a trisp closure map: The simplex $(b_{1}<b_{2}<b_{3}<b_{4})$
cannot be extended by $r_{3}$, and the situation is not remedied
by considering the dual poset instead. On the other hand, Figure~\ref{fig:counterex}b
gives an example of a trisp closure map that is not induced by any
order-preserving map.

If we shift our point of view from posets to acyclic categories, as
arises naturally when asking about a trisp closure map on $\Delta(P/G)$,
we loose this tool. But we can still ask for any AC-map $\phi:C\rightarrow C$
whether the induced map on the vertices of $\Delta(C)$ is a trisp
closure map for some choice of $B$. Necessary conditions are:
\begin{condition}
\label{con:ac-map to tr cl map}$B$ and $\phi(B)$ are disjoint,
and $|\mathcal{M}(x,\phi(x))|+|\mathcal{M}(\phi(x),x)|=1$ for all
$x\in B$.
\end{condition}

\section{\label{sec:Cl-maps-on-quot-trisps}Closure maps on $\Delta(C)/G$}

Consider first some arbitrary trisp $T$ with a $G$-action, where
$G$ is a finite group. According to Kozlov \cite[Section 14.1.2]{Ko08a},
the following condition ensures that $T/G$ is again a regular trisp.
\begin{condition*}
[R]For any $g\in G$ and any simplex $\sigma$ of $T$, the simplex
$g(\sigma)\cap\sigma$ is fixed pointwise by $g$.
\end{condition*}
It turns out that this condition ensures that the quotient map of
a $G$-equivariant trisp closure map on $T$ is a trisp closure map
on $T/G$.
\begin{prop}
\label{pro:trispcl=00003D>trispGcl}Let $T$ be a finite trisp with
a $G$-action fulfilling Condition~R. Let $\phi:B\rightarrow R$
be a $G$-equivariant trisp closure map on $T$, where $B$ and $R$
are closed under $G$ (that is, $G\cdot B\cap G\cdot R=\emptyset$).
Then the quotient map $\phi_{G}$ is a trisp closure map on $T/G$.\end{prop}
\begin{proof}
Let $G\sigma$ be a simplex of $T/G$, with a minimal blue vertex
$Gb$. Choose the representative $b$ such that it is a vertex of
$\sigma$; since the group action respects trisp order, $b$ is the
minimal blue vertex of $\sigma$ as well. There is only something
to prove if $\phi_{G}(Gb)=G\phi(b)$ is not a vertex of $G\sigma$.
Then $\phi(b)$ is not a vertex of $\sigma$ and there exists a unique
extension simplex $\tau$ of $\sigma$ by $\phi(b)$; that is, $\phi(b)$
is the $j$-th vertex of $\tau$ and $\partial_{j}\tau=\sigma$. Then
$\partial_{j}(G\tau)=G\sigma$ and $G\phi(b)$ is a vertex of $G\tau$,
so $G\tau$ extends $G\sigma$.

Assume that there is another extension $G\tau'$ of $G\sigma$. Choose
the representative $\tau'$ such that $\partial_{j}\tau'=\sigma$
for some $j$, then $\tau'$ contains the vertices $b$ and $g\phi(b)$,
where $g\in G$ such that $g\phi(b)$ is the representative of $G\phi(b)$
in $\tau'$. Thus there is a simplex $\rho$ in $T$ with vertex set
$\{b,g\phi(b)\}$. Assume that $g\phi(b)\neq\phi(b)$, since otherwise
$\tau'=\tau$. As $\phi(b)$ is a red vertex and thus $g\phi(b)$
is red as well, there must exist an extension simplex of $\rho$,
with vertices $b,\phi(b),g\phi(b)$, and a boundary simplex $\rho'$
with vertices $\phi(b)$ and $g\phi(b)$. Then $g\phi(b)\in g\rho'\cap\rho'$,
so by Condition~S we have $g\phi(b)=g^{2}\phi(b)$ and thus $\phi(b)=g\phi(b)$.
Hence $\tau=\tau'$ and thus the extension simplex of $G\sigma$ is
unique in $T/G$.
\end{proof}
As shown in \cite[Prop. 14.4]{Ko08a}, Condition~R is automatically
fulfilled for $T=\Delta(C)$ for some acyclic category $C$, where
the group action on $T$ is induced by a group action on $C$. So
we obtain the following corollary.
\begin{cor}
\label{cor:ac-map->G-tr-cl-map}Let $C$ be a finite acyclic category
with $G$-action. Let $\phi:C\rightarrow C$ be a $G$-equivariant
AC-map that induces a trisp closure map $\bar{\phi}:B\rightarrow R$
on $\Delta(C)$, where $B$ and $R$ are closed under $G$ (that is,
$G\cdot B\cap G\cdot R=\emptyset$). Then the quotient map $\bar{\phi}_{G}$
is a trisp closure map on $\Delta(C)/G$.
\end{cor}
Now we consider the opposite direction: Given a trisp $T$ with a
group action $G$ fulfilling Condition~R (as we want to have regularity
of $T/G$ ensured), and a trisp closure map $\psi:\bar{B}\rightarrow\bar{R}$
on $T/G$, we wish to obtain a trisp closure map $\phi:B\rightarrow R$
on $T$ such that $\phi_{G}=\psi$. There are no choices for $B=\{b\in X|X\in\bar{B}\}$
and $R=\{r\in X|X\in\bar{R}\}$.
\begin{condition*}
[C]For each $b\in B$ there is a unique vertex $r_{b}\in\psi(Gb)$
such that there exists a simplex with vertex set $\{b,r_{b}\}$ and
this simplex is unique as well.
\end{condition*}
In other words, if Condition~C is fulfilled, then $\phi$ is well-defined
by setting $\phi(b)=r_{b}$. 
\begin{prop}
\label{pro:ConC is necessary}Condition~C is necessary for the existence
of a trisp closure map $\phi:B\rightarrow R$ on $T$ with $\phi_{G}=\psi$.\end{prop}
\begin{proof}
For each $b\in B$ the image $\phi(b)=:r$ must be chosen from $\psi(Gb)$.
As $\psi$ is a closure map, there exists exactly one simplex $G\sigma$
in $T/G$ with vertex set $\{Gb,\psi(Gb)\}$, and a representative
$\sigma$ can be chosen with $b$ as a vertex. The other vertex of
$\sigma$ is then $r$. If there exists another candidate $gr\in\psi(Gb)$
such that there is a simplex $\sigma'$ with vertex set $\{b,gr\}$,
then an extension $\tau$ of $\sigma'$ must exist, with vertex set
$\{b,gr,\phi(b)=r\}$. One boundary simplex $\tau'$ of $\tau$ has
vertex set $\{gr,r\}$, which by the route of $g\tau'\cap\tau'\ni gr$
and Condition~R implies $gr=r$.

If there are two different simplices $\sigma,\sigma'$ with vertex
set $\{b,r\}$, then both extend the simplex $b$, in contradiction
to $\phi$ being a trisp closure map.
\end{proof}
In general, Condition~C is not sufficient to obtain $\phi$. Consider
for example the following regular trisp: 0-simplices are $b,x,r$;
1-simplices are $(b,x),(b,r),(x,r)$ and two 2-simplices $\sigma$
and $\tau$, both having all the 1-simplices as boundaries (a filled
triangle with double filling). $\mathbb{Z}_{2}$ acts on this trisp:
the nonidentity element interchanges $\sigma$ and $\tau$, leaving
all other simplices fixed. This action is a trisp action fulfilling
Condition~C. The quotient trisp is just a filled triangle, and $\psi:\{\mathbb{Z}_{2}b\}\rightarrow\{\mathbb{Z}_{2}x,\mathbb{Z}_{2}r\}$,
mapping $\mathbb{Z}_{2}b$ to $\mathbb{Z}_{2}r$, is a $\mathbb{Z}_{2}$-equivariant
trisp closure map that is not the quotient of any trisp closure map
in the original trisp.

On the other hand, if $T$ is an abstract simplicial complex (meaning
here that different simplices have different vertex sets), then Condition~C
is indeed sufficient.
\begin{prop}
Let $T$ be an abstract simplicial complex with a group $G$ acting
on $T$ under Condition~R; let $\psi:\bar{B}\rightarrow\bar{R}$
be a trisp closure map on $T/G$. Then the following are equivalent:

(1) $\psi$ fulfills Condition~C

(2) $\phi:B\rightarrow R$ is a trisp closure map on $T$ with $\phi_{G}=\psi$,
where $B=\{b\in X|X\in\bar{B}\}$ and $R=\{r\in X|X\in\bar{R}\}$
and $\phi(b)=r_{b}$, with $r_{b}$ as in Condition~C.\end{prop}
\begin{proof}
Condition~C is necessary by Proposition~\ref{pro:ConC is necessary}.
To show that it is sufficient, consider a simplex $\sigma$ of $T$,
with minimal blue vertex $b$ and $\phi(b)\notin\sigma$. Then $Gb$
is the minimal blue vertex of $G\sigma$ and there exists an extension
$G\tau$ of $G\sigma$. So $j$ exists with $G\sigma=\partial_{j}(G\tau)=G(\partial_{j}\tau)$.
Choose a representative $\tau$ such that $\partial_{j}\tau=\sigma$,
which is possible by the regularity of the $G$-action on $T$ implied
by Condition~R. Then the $j$-th vertex $v_{j}$ of $\tau$ is in
$\psi(Gb)$, and there is a subsimplex $\{b,v_{j}\}$ of $\tau$,
so by Condition~C we have $v_{j}=\phi(b)$. Thus $\tau$ is the unique
extension simplex of $\sigma$.
\end{proof}
In particular, the nerve of a poset $P$ is always an abstract simplicial
complex, and a group $G$ acting on $P$ induces an action on $\Delta(P)$
that always fulfills Condition~R.

\section{\label{sec:cl-maps-on-quot-poset}closure maps on $\Delta(P/G)$}

As $P/G$ is usually not a poset, we cannot hope that the quotient
of a closure operator is again a closure operator. But it turns out
that the induced map on $\Delta(P/G)$ is still a trisp closure map.

All results in this section hold for ascending closure operators as
well, using the {}``maximal'' version of trisp closure maps.
\begin{lem}
\label{lem:phi(morph)}Let $\phi:P\rightarrow P$ be a $G$-equivariant
descending closure operator. Then $[b<a]=[c<a]$ in $P/G$ implies
that $[\phi(b)<\phi(a)]=[\phi(c)<\phi(a)]$ in $P/G$ and also in
$\phi(P)/G$.\end{lem}
\begin{proof}
Long version: Writing this out as in Section~\ref{sub:The-quotient-C/G},
this means that there exist $n,t_{1},\ldots,t_{n}\in\mathcal{N},a_{ij}^{+},a_{ij}^{-}\in P,g_{ij}\in G$
such that $b=a_{11}^{+}<a_{12}^{+}<\ldots<a_{1t_{1}}^{+}=a$, $a_{i1}^{-}<a_{i2}^{-}<\ldots<a_{it_{i-1}}^{-}$,
$a_{i1}^{+}<a_{i2}^{+}<\ldots<a_{it_{i}}^{+}$, $c=a_{n1}^{-}<a_{n2}^{-}<\ldots<a_{nt_{n-1}}^{-}=a$,
with $a_{i1}^{-}=a_{i1}^{+}$, $a_{it_{i-1}}^{-}=a_{it_{i}}^{+}$,
$g_{ij}a_{ij}^{+}=a_{i+1,j}^{-}$, $g_{ij}a_{i,j+1}^{+}=a_{i+1,j+1}^{-}$.

An example of the situation is shown in Figure~\ref{fig:n=00003D3-ex}.

\begin{figure}
\includegraphics{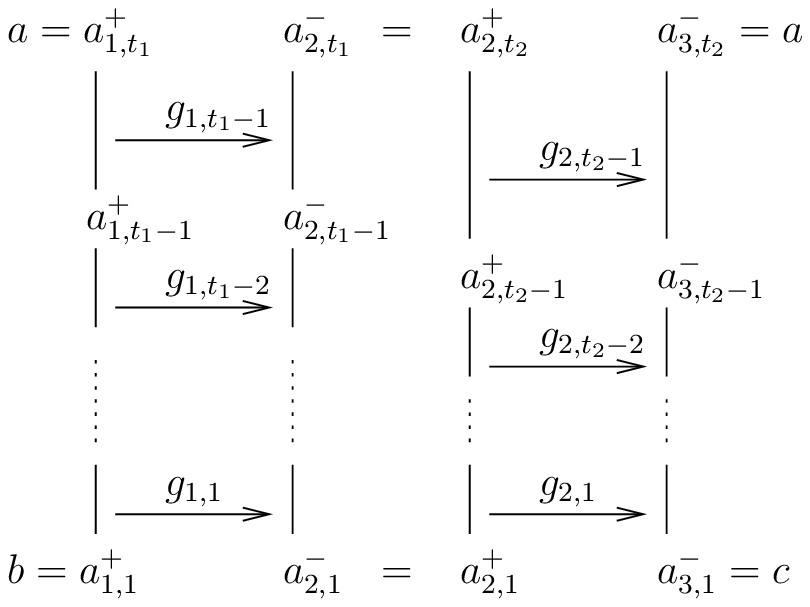}

\caption{\label{fig:n=00003D3-ex}Example situation with $n=3$}

\end{figure}

We will prove the statement by induction on $n$. 

If $n=2$, then by applying $\phi$ to the chains, we obtain $\phi(b)=\phi(a_{11}^{+})\leq\phi(a_{12}^{+})\leq\ldots\leq\phi(a_{1t_{1}}^{+})=\phi(a)$,
$\phi(c)=\phi(a_{21}^{-})\leq\phi(a_{22}^{-})\leq\ldots\leq\phi(a_{2t_{1}}^{-})=\phi(a)$.
In each of these, we might get subsequences with equality $\phi(a_{jk})=\phi(a_{j,k+1})=\ldots=\phi(a_{j,k+l})$.
Reduce the index set $\{1,2,\ldots,t_{1}\}$ by keeping only the first
index in each of these sequences; this yields the same result for
both chain images because the $G$-action is horizontal. To simplify
notation, denote the new index set with $\{1,2,\ldots,t_{1}\}$ as
well. Thus we get chains $\phi(b)=\phi(a_{11}^{+})<\phi(a_{12}^{+})<\ldots<\phi(a_{1t_{1}}^{+})=\phi(a)$,
$\phi(c)=\phi(a_{21}^{-})<\phi(a_{22}^{-})<\ldots<\phi(a_{2t_{1}}^{-})=\phi(a)$,
with $g_{ij}\phi(a_{ij}^{+})=\phi(g_{ij}a_{ij}^{+})=\phi(a_{i+1,j}^{-})$,
$g_{ij}\phi(a_{i,j+1}^{+})=\phi(a_{i+1,j+1}^{-})$. We conclude that
$[\phi(b)<\phi(a)]=[\phi(c)<\phi(a)]$.

Assume that $[\phi(b)<\phi(a)]=[\phi(a_{n-1,1}^{-})<\phi(a_{n-1,t_{n-2}}^{-})]$.
Proceeding as above using the chains $a_{n-1,1}^{+}<a<\ldots<a_{n-1,t_{n-1}}^{+}$,
$c=a_{n1}^{-}<a_{n2}^{-}<\ldots<a_{nt_{n-1}}^{-}=a$ yields $[\phi(c)<\phi(a)]=[\phi(a_{n-1,1}^{+})<\phi(a_{n-1,t_{n-1}}^{+})]=[\phi(a_{n-1,1}^{-})<\phi(a_{n-1,t_{n-2}}^{-})]=[\phi(b)<\phi(a)]$.\end{proof}
\begin{prop}
\label{pro:desc-image-subtrisp}Let $\phi:P\rightarrow P$ be a $G$-equivariant
descending closure operator. Then $\Delta_{R}(P/G)$, the subtrisp
of $\Delta(P/G)$ induced by the vertex set $R=G\cdot\phi(P)$, equals
$\Delta(\phi(P)/G)$.\end{prop}
\begin{proof}
The vertex sets of both trisps coincide by definition. By Lemma~\ref{lem:phi(morph)},
the 1-skeleta coincide as well, and because both trisps are flag complexes,
this is sufficient to guarantee equality.\end{proof}
\begin{prop}
\label{pro:desc-cl-to-quot-poset}Let $\phi:P\rightarrow P$ be a
$G$-equivariant descending closure operator. Then the induced map
$\phi_{G}$ is a trisp closure map on $\Delta(P/G)$.\end{prop}
\begin{proof}
Let $R=\phi_{G}(P/G)$ be the red vertices of $P/G$, $B=(P/G)\backslash R$
be the blue vertices. Consider some simplex $\sigma$ containing a
minimal blue vertex $[b]$ and not containing $\phi_{G}([b])$. There
is a $\lambda$-preimage $G(a_{1}<a_{2}<\ldots<a_{r})$ of $\sigma$,
which can be extended by $[\phi(b)]=\phi_{G}([b])$, since the induced
map on $\Delta(P)/G$ is a trisp closure map by Remark~\ref{rem:desc-cl-map->tr-cl-map}
and Proposition~\ref{pro:trispcl=00003D>trispGcl}. Taking the $\lambda$-image
of this extension provides an extension of $\sigma$ in $\Delta(P/G)$. 

The minimal blue vertex in $\sigma$ is $[b]=[a_{i}]$ and $\phi$
is descending, so there can be no subsimplex $([\phi(b)]\stackrel{[x]}{\rightarrow}[a_{i-1}]\stackrel{[m_{i-1}]}{\rightarrow}[b])$
of an extension of $\sigma$. Otherwise we could obtain a preimage
$g\phi(b)<a_{i-1}<b$ by properly choosing representatives $b,a_{i+1},g\phi(b)$.
Applying $\phi$ leads to $g\phi(b)\leq\phi(a_{i+1})=a_{i+1}\leq\phi(b)$,
thus $\phi(b)=a_{i+1}=g\phi(b)$ . Therefore any extension of $\sigma$
must be of the form $([a_{1}]\stackrel{[m_{1}]}{\rightarrow}[a_{2}]\stackrel{[m_{2}]}{\rightarrow}\ldots\stackrel{[m_{i-2}]}{\rightarrow}[a_{i-1}]\stackrel{[y]}{\rightarrow}[\phi(b)]\stackrel{[x]}{\rightarrow}[b]=[a_{i}]\stackrel{[m_{i}]}{\rightarrow}\ldots\stackrel{[m_{r-1}]}{\rightarrow}[a_{r}])$.

The only remaining question is whether there are different extensions
to $\sigma$. These can vary only in their choice of $[x]$ and $[y]$,
under the condition that $[x\circ y]=[m_{i-1}]$. So we only need
to consider the situation where $[r]\stackrel{[y]}{\rightarrow}[\phi(b)]\stackrel{[x]}{\rightarrow}[b]$
with $[r]\in R$, and $[r]\stackrel{[y']}{\rightarrow}[\phi(b)]\stackrel{[x']}{\rightarrow}[b]$
with $[x'\circ y']=[x\circ y]$ and prove that then $[x]=[x']$ and
$[y]=[y']$.

(1) Let $g\in G$ such that $[x]=[g\phi(b)<b]$. Since $g\phi(b)=\phi(g\phi(b))\leq\phi(b)$,
we see that $[x]=[\phi(b)<b]$; by the same argument $[x']=[\phi(b)<b]$
holds.

(2) Choose a representative $r$ such that $[y]=[r<\phi(b)]$, let
$g\in G$ such that $[y']=[gr<\phi(b)]$. Thus $[x\circ y]=[r<b]$,
which by our assumptions equals $[x'\circ y']=[gr<b]$. By Lemma~\ref{lem:phi(morph)},
we have $[y]=[\phi(r)<\phi(b)]=[\phi(gr)<\phi(b)]=[y']$.
\end{proof}

\section{\label{sec:Applications}Applications}

Vassiliev introduced in his work on knot invariants (\cite{Va93})
the complexes of disconnected graphs. The vertex set of such a complex
$DG_{n}$ consists of all 2-element subsets of $\{1,\ldots,n\}$,
indexing all possible edges in a graph on $n$ vertices. The simplices
of $DG_{n}$ are exactly the edge sets of disconnected graphs on $n$
vertices. $DG_{n}$ carries a $\mathcal{S}_{n}$-action induced by
the action on the graph vertices, though this action does not fulfill
Condition~R. But there is an induced $\mathcal{S}_{n}$-action on
the face poset $\mathcal{F}(DG_{n})$, which in turn lets us explore
the trisp $\Delta(\overline{\mathcal{F}(DG_{n})})/\mathcal{S}_{n}=\bd(DG_{n})/\mathcal{S}_{n}$
with our tools, which simplifies the first analysis of this trisp
by Kozlov in \cite{Ko08b}.

Let $\overline{\Pi_{n}}$ be the poset of partitions of $\{1,\ldots,n\}$
ordered by refinement except $1|2|\ldots|n$ and $12\ldots n$, let
$\phi:\overline{\mathcal{F}(DG_{n})}\rightarrow\overline{\mathcal{F}(DG_{n})}$
be the map taking each graph $G$ to its transitive closure, that
is the direct sum of the complete graphs on each of the components
of $G$. So $\phi(G)$ can be understood as a partition of $\{1,\ldots,n\}$
and $\phi(\overline{\mathcal{F}(DG_{n}})$ is isomorphic to $\overline{\Pi_{n}}$.
\begin{cor}
The trisp $\bd(DG_{n})/\mathcal{S}_{n}$ is collapsible.\end{cor}
\begin{proof}
Note that $\phi$ is a $G$-equivariant ascending closure operator,
so by Remark~\ref{rem:desc-cl-map->tr-cl-map}, the prerequisites
for Corollary~\ref{cor:ac-map->G-tr-cl-map} are fulfilled. Hence
$\bd(DG_{n})/\mathcal{S}_{n}$ collapses onto $\Delta(\overline{\Pi_{n}})/\mathcal{S}_{n}$,
the collapsibility of which has been shown in \cite{Ko00}.
\end{proof}
Our results allow us to tackle the nerve of $\overline{\mathcal{F}(DG_{n})}/\mathcal{S}_{n}$
using the same map $\phi$:
\begin{cor}
The trisp $\Delta(\overline{\mathcal{F}(DG_{n})}/\mathcal{S}_{n})$
is collapsible.\end{cor}
\begin{proof}
By Propositions~\ref{pro:desc-cl-to-quot-poset} and \ref{pro:desc-image-subtrisp},
$\Delta(\overline{\mathcal{F}(DG_{n})}/\mathcal{S}_{n})$ collapses
onto $\Delta(\overline{\Pi_{n}}/\mathcal{S}_{n})$. The objects of
the category $\overline{\Pi_{n}}/\mathcal{S}_{n}$ can be indexed
by the nontrivial number partitions of $n$, between which $2+1+\ldots+1$
is minimal. We will show that this is in fact a terminal object, hence
$\overline{\Pi_{n}}/\mathcal{S}_{n}$ is collapsible.

Let $\pi=\pi_{1}|\pi_{2}|\ldots\in\overline{\Pi_{n}}$, denote by
$\overline{ab}$ the partition with one set $\{a,b\}$, rest singletons.
Fix $a,b\in\pi_{1}$, let $\overline{cd}$ be another refinement of
$\pi$. If $c,d$ are in $\pi_{1}$ as well, then $g:=(ac)(bd)$ stabilizes
$\pi$, thus $[\overline{ab}<\pi]=[\overline{cd}<\pi]$. If $c,d$
are in some other set $\pi_{i}$, then we decompose $\overline{ab}<\pi$
and $\overline{cd}<\pi$ into $\overline{ab}<\overline{ab|cd}<\pi$
and $\overline{cd}<\overline{ab|cd}<\pi$, and note that $(ac)(bd),id$
map the former to the latter. Hence we obtain $[\overline{ab}<\pi]=[\overline{cd}<\pi]$,
that is, for each $[\pi]\in\overline{\Pi_{n}}/\mathcal{S}_{n}$ there
is only a single morphism with source $[\pi]$ and target $2+1+\ldots+1$.\end{proof}
\begin{acknowledgement*}
I would like to thank Prof. Feichtner-Kozlov for interesting me in
this research.
\end{acknowledgement*}

\end{document}